\documentclass[12pt]{article} 
\usepackage{amsmath,amssymb}
\usepackage{amscd}
\usepackage[all]{xy} \newtheorem{thm}{Theorem}
\newtheorem{df}{Def\/inition} \newtheorem{prop}{Proposition}
\newtheorem{cor}{Corollary} \newtheorem{lemma}{Lemma}

\DeclareMathOperator{\ord}{ord}
\DeclareMathOperator{\id}{id} 
 \DeclareMathOperator{\Hom}{Hom}
\DeclareMathOperator{\Spec}{Spec} \DeclareMathOperator{\Ker}{Ker}
\DeclareMathOperator{\Coker}{Coker} 
 
\DeclareMathOperator{\Ext}{Ext}

\def\endproof{$\hfill \square$}

\begin{document}
\title{Boundedness results for finite flat group schemes over discrete
  valuation rings of mixed characteristic} 
\author{Adrian Vasiu and
  Thomas Zink}
\maketitle

\centerline{to appear in J. Number Theory}

\medskip\noindent {\bf Abstract.} Let $p$ be a prime. Let $V$ be a
discrete valuation ring of mixed characteristic $(0,p)$ and index of
ramification $e$. Let $f: G \rightarrow H$ be a homomorphism of finite flat
commutative group schemes of $p$ power order over $V$ whose generic
fiber is an isomorphism. We provide a new proof of a result of Bondarko and Liu that bounds the kernel and the cokernel of the special fiber of $f$ in terms of $e$.  For $e < p-1$ this reproves a result of Raynaud. Our bounds are sharper that the ones of Liu, are almost as sharp as the ones of Bondarko, and involve a very simple and short method. As an application we obtain a new proof of an extension theorem for homomorphisms of truncated Barsotti--Tate groups which strengthens Tate's extension theorem for homomorphisms of $p$-divisible groups. 



\section{Introduction}
Let $p$ be a rational prime. Let $V$ be a discrete valuation ring of
mixed characteristic $(0,p)$. Let $K$ be the field of fractions of
$V$. Let $k$ be the residue field of $V$. Let $e$ be the index of
ramification of $V$. Let $G$ and $H$ be two finite flat commutative
group schemes of $p$ power order over $V$. For $n\in
\mathbb{N}^{\ast}$, let $G\lceil p^n\rceil$ be the schematic closure
of $G_K[p^n]$ in $G$. The goal of the paper is to provide new proofs using {\it Breuil modules} of the following theorem and of several applications of it.

\begin{thm}\label{T1} There exists a non-negative integer $s$ that
  depends only on $V$ and that has the following property.

  For each homomorphism $f:G\rightarrow H$ whose generic fiber $f_K:G_K\rightarrow H_K$
  is an isomorphism, there exists a homomorphism $f':H\rightarrow G$ such
  that $f'\circ f=p^s\id_G$ and $f\circ f'=p^s\id_H$ and
  therefore the special fiber homomorphism $f_k:G_k\rightarrow H_k$ has a
  kernel and a cokernel annihilated by $p^s$. If moreover $H$ is a
  truncated Barsotti--Tate group of level $n > s$, then the natural
  homomorphism $f\lceil p^{n-s}\rceil:G\lceil p^{n-s}\rceil\rightarrow
  H[p^{n-s}]$ is an isomorphism.
\end{thm}

An equivalent form of the first part of Theorem \ref{T1} (see the first part of Corollary 2 below) was first obtained in \cite{Bon} using {\it Cartier--Dieudonn\'e modules}, in \cite{Liu1} using {\it deformation theories} and {\it Galois-descent}, and in \cite{Liu2} using {\it Breuil--Kisin modules}. The number $s$ admits computable upper bounds in terms only of $e$. For instance, we have $s\leq (\log_p e + \ord_p e +2)(\ord_p e +2)$ (cf. Examples 2 and 4). This upper bound of $s$ is weaker than the one obtained by Bondarko $s\leq\lfloor \log_p({{pe}\over {p-1}})\rfloor$ but it is much stronger than the upper bounds one can get based on the proofs of \cite{Liu1}, Theorem 1.0.5 and \cite{Liu2}, Theorem 2.4.2. If $p$ does not divide $e$, then we regain Bondarko's upper bound (cf. Examples 1 and 2). If $e\le p-2$, then $s=0$ (cf. Example 1) and thus we also regain the following classical theorem of Raynaud.

\begin{cor}\label{C1} 
  We assume that $e\le p-2$ (thus $p$ is odd). Then each finite flat
  commutative group scheme of $p$ power order over $K$ extends in at
  most one way to a finite flat commutative group scheme over $V$.
\end{cor}

Corollary \ref{C1} was first proved in \cite{R}, Theorem 3.3.3 or
Corollary 3.3.6 and more recently in \cite{Bon}, \cite{Liu1}, \cite{Liu2}, and \cite{VZ2}, Proposition 15. The
first part of the next result is an equivalent form of the first part
of Theorem \ref{T1}. 

\begin{cor}\label{C2}
Let $h: G_K \rightarrow  H_K$ be a homomorphism over $K$. Then $p^s h$
extends to a homomorphism $G \rightarrow H$ (i.e., the cokernel of the
natural monomorphism $\Hom(G,H)\hookrightarrow \Hom(G_K,H_K)$ is
annihilated by $p^s$). Thus the natural homomorphism
$\Ext^1(H,G)\rightarrow \Ext^1(H_K,G_K)$ has a kernel annihilated by
$p^s$.   
\end{cor}

The first (resp. the second) part of Corollary \ref{C2} was first obtained in \cite{Bon}, Theorem A or B, in \cite{Liu1}, Theorem 1.0.5, and in \cite{Liu2}, Theorem 2.4.2 (resp.  in \cite{Bon}, Theorem D). 

The next two results are a mixed characteristic geometric
analogue of the homomorphism form \cite{V}, Theorem 5.1.1 of the {\it
  crystalline boundedness principle} presented in \cite{V}, Theorem
1.2. 

\begin{cor}\label{C3}
  We assume that $G$ and $H$ are truncated Barsotti--Tate groups of
  level $n>s$. Let $h:G_K\rightarrow H_K$ be a homomorphism. Then the
  restriction homomorphism $h[p^{n-s}]:G_K[p^{n-s}]\rightarrow H_K[p^{n-s}]$
  extends to a homomorphism $G[p^{n-s}]\rightarrow H[p^{n-s}]$.
\end{cor}

\begin{cor}\label{C4}
  We assume that $n>2s$. Let $H$ be a truncated Barsotti--Tate group
  of level $n$ over $V$. Let $G$ be such that we have an isomorphism
  $h:G_K\rightarrow H_K$. Then the quotient group scheme $G\lceil
  p^{n-s}\rceil/G\lceil p^{s}\rceil$ is isomorphic to $H[p^{n-2s}]$
  and thus it is a truncated Barsotti--Tate group of level $n-2s$.
\end{cor}

A weaker form of Corollary \ref{C3} was first obtained in \cite{Bon}, Theorem E. 

Let $Y$ be a normal noetherian integral scheme with field of functions
$L$ of characteristic zero. A classical theorem of Tate (\cite{T},
Theorem 4) says that for every two $p$-divisible groups $D$ and $F$
over $Y$, each homomorphism $D_L \rightarrow F_L$ extends uniquely to
a homomorphism $D\rightarrow F$. From Corollary \ref{C3} we obtain the
following sharper version of this theorem.

\begin{thm}\label{T2}
  Let $Y$ and $L$ be as above. Then there exists a non-negative
  integer $s_Y$ which has the following property.

  Let $\mathcal G$ and $\mathcal H$ be truncated Barsotti--Tate groups
  over $Y$ of level $n > s_Y$ and of order a power of the prime
  $p$. Let $h: \mathcal G_L \rightarrow \mathcal H_L$ be a
  homomorphism. Then there exists a unique homomorphism $g: \mathcal
  G[p^{n-s_Y}] \rightarrow \mathcal H[p^{n-s_Y}]$ that induces
  $h[p^{n-s_Y}]$ over $L$.
\end{thm}

Section 2 recalls the classification of
finite flat commutative group schemes of
$p$ power order over $V$ in terms of Breuil modules that holds
provided $k$ is perfect and $V$ is complete. This classification was
conjectured by Breuil (see \cite{Br}), was first proved in \cite{K1} and \cite{K2} for $p>2$ and for the connected case with $p=2$, was generalized (using a covariant language) in
\cite{VZ1}, \cite{Lau1}, and \cite{Lau2} for all primes $p$, and was recently proved for $p=2$ in \cite{Lau2}, \cite{Liu3}, and \cite{Kim}. 
 As in \cite{Liu2}, our method to prove the above results relies as well on Breuil modules but it performs a refined combinatorial analysis of {\it Eisenstein polynomials} associated to uniformizers of $V$ which allows us to get better upper bounds than the ones of Liu and simpler proofs (see Lemmas 1 to 4 and Proposition 2) that could be generalized to the crystalline context of \cite{Liu2}, Theorem 2.4.2. 
In Section 3 we provide formulas
for $s$ as well as explicit upper bounds of it. In Section 4 we prove
Theorem \ref{T1}. Corollaries 1 to 4 and Theorem \ref{T2} are
proved in Section 5. In Section 6 we present extra applications to
different heights associated to $G$.

\section{Breuil modules}\label{2}

Let $V \hookrightarrow V'$ be an extension of discrete valuation rings
of mixed characteristic $(0,p)$ such that the index of ramification of $V'$ is also $e$. Theorem
\ref{T1} holds for $V$ if it holds for $V'$. There exists an extension
$V'$ which is complete and has a perfect residue field. If we find an
upper bound of $s$ only in terms of $e$ which holds for each complete
$V'$ with perfect residue field, then this upper bound of $s$ is also
good for $V$.

Thus from now on we assume that $V$ is complete and that $k$ is a
perfect field. Let $W(k)$ be the ring of Witt vectors with
coefficients in $k$. We will view $V$ as a $W(k)$-algebra which is a
free $W(k)$-module of rank $e$. Let $\ord_p:W(k)\rightarrow
\mathbb{N}\cup\{\infty\}$ be the $p$-adic valuation normalized by the
conditions that $\ord_p(p)=1$ and $\ord_p(0)=\infty$. Let $u$ be a
variable and let
\begin{displaymath}
  \mathfrak{S}:=W(k)[[u]].
\end{displaymath}
We extend the Frobenius endomorphism $\sigma$ of $W(k)$ to
$\mathfrak{S}$ by the rule $\sigma(u)=u^p$. For
$n\in\mathbb{N}^{\ast}$ let
$\mathfrak{S}_n:=\mathfrak{S}/p^n\mathfrak{S}$. If $M$ is a
$\mathfrak{S}$-module let
\begin{displaymath}
  M^{(\sigma)} := \mathfrak{S} \otimes_{\sigma, \mathfrak{S}} M. 
\end{displaymath} 

Let $\pi$ be a uniformizer of $V$. Let
\begin{displaymath}
  E=E(u)=u^e+a_{e-1}u^{e-1}+\dots+ a_1 u+a_0\in W(k)[u]
\end{displaymath}
be the unique Eisenstein polynomial in $u$ which has coefficients in
$W(k)$ and which has $\pi$ as a root. For $i\in\{0,\ldots,e-1\}$ we
have $a_i\in pW(k)$; moreover $a_0$ is $p$ times a unit of $W(k)$. We
have a $W(k)$-epimorphism
$$q_{\pi}:\mathfrak{S}\twoheadrightarrow V$$
that maps $u$ to $\pi$.

\begin{df}\label{D1}
  By a (contravariant) Breuil window relative to
  $q_{\pi}:\mathfrak{S}\twoheadrightarrow V$ we mean a pair
  $(Q,\phi)$, where $Q$ is a free $\mathfrak{S}$-module of finite rank
  and where $\phi: Q^{(\sigma)}\rightarrow Q$ is a $\mathfrak{S}$-linear map
  (Frobenius map) whose cokernel is annihilated by $E$.  By a
  (contravariant) Breuil module relative to
  $q_{\pi}:\mathfrak{S}\twoheadrightarrow V$ we mean a pair
  $(M,\varphi)$, where $M$ is a finitely generated $\mathfrak{S}$-module annihilated by a
  power of $p$ and of projective dimension at most one and where
  $\varphi: M^{(\sigma)}\rightarrow M$ is a $\mathfrak{S}$-linear map
  (Frobenius map) whose cokernel is annihilated by $E$. 
\end{df}

\begin{df}\label{D2}
  Let $\mathcal B$ be the category of
  Breuil modules relative to
  $q_{\pi}:\mathfrak{S}\twoheadrightarrow V$. Let $\mathcal B_1$
  be the full subcategory of
  $\mathcal B$ whose objects are
  Breuil modules $(M,\varphi)$ relative to
  $q_{\pi}:\mathfrak{S}\twoheadrightarrow V$ with $M$ annihilated by
  $p$. Let $\mathcal F$ be the category of
  finite flat commutative group schemes of $p$ power order over $V$. Let
  $\mathcal F_1$ be the full subcategory of
  $\mathcal F$ whose objects are finite
  flat commutative group schemes over $V$ annihilated by $p$.
\end{df}

We explain shortly the relation to covariant Breuil windows and
modules relative to
  $q_{\pi}:\mathfrak{S}\twoheadrightarrow V$ as introduced in \cite{VZ1}, Definitions 1 and 2. Let $L$ be a free
$\mathfrak{S}$-module of finite rank. We have a canonical $\mathfrak{S}$-linear isomorphism 
\begin{displaymath}
   \mathfrak{S} \otimes_{\sigma, \mathfrak{S}} \Hom_{\mathfrak{S}}(L,
   \mathfrak{S}) \rightarrow \Hom_{\mathfrak{S}}(L^{(\sigma)},
   \mathfrak{S})
\end{displaymath}
that maps $s \otimes \alpha$ to the $\sigma$-linear map $L
\rightarrow  \mathfrak{S}$, $l \mapsto s \sigma(\alpha(l))$. The last
$\sigma$-linear map may be also regarded as an element of $\Hom_{\mathfrak{S}}(L^{(\sigma)},
   \mathfrak{S})$. As $\mathfrak{S}
   \overset{\sigma}{\longrightarrow} \mathfrak{S}$ is flat, we obtain
   more generally that for each $\mathfrak{S}$-module $M$ of finite type
   and for every integer $i\geq 0$ we have a canonical $\mathfrak{S}$-linear isomorphism 
\begin{displaymath}
   \Ext_{\mathfrak{S}}^i(M, \mathfrak{S})^{(\sigma)} \cong
   \Ext_{\mathfrak{S}}^i(M^{(\sigma)}, \mathfrak{S}).
\end{displaymath} 

Let $(Q, \phi)$ be a contravariant Breuil window relative to
  $q_{\pi}:\mathfrak{S}\twoheadrightarrow V$. Applying the functor
$\Hom_{\mathfrak{S}}(\,-\, , \mathfrak{S})$ to the $\mathfrak{S}$-linear map $\phi$
we obtain a homomorphism
\begin{displaymath}
  \psi: \Hom_{\mathfrak{S}}(Q, \mathfrak{S}) \rightarrow
  \Hom_{\mathfrak{S}}(Q^{(\sigma)}, \mathfrak{S}) \cong \Hom_{\mathfrak{S}}(Q,
  \mathfrak{S})^{(\sigma)}. 
\end{displaymath} 
The pair $(\Hom_{\mathfrak{S}}(Q, \mathfrak{S}),\psi)$ is a covariant Breuil window relative to
  $q_{\pi}:\mathfrak{S}\twoheadrightarrow V$. The rule on objects $(Q,\phi)\mapsto (\Hom_{\mathfrak{S}}(Q, \mathfrak{S}),\psi)$ induces an
antiequivalence of the category of contravariant Breuil windows relative to
  $q_{\pi}:\mathfrak{S}\twoheadrightarrow V$ with
the category of covariant Breuil windows relative to
  $q_{\pi}:\mathfrak{S}\twoheadrightarrow V$. 

In the same way the functor $\Ext_{\mathfrak{S}}^1(\, - \, , \mathfrak{S})$
induces naturally an antiequivalence of the category of contravariant Breuil
modules relative to
  $q_{\pi}:\mathfrak{S}\twoheadrightarrow V$ with the category of covariant Breuil modules relative to
  $q_{\pi}:\mathfrak{S}\twoheadrightarrow V$.

If $M$ is annihilated by $p$, then it is easy to see that $M$ is a
free $\mathfrak{S}_1$-module (cf. \cite{VZ2}, Section 2, p. 578); its
rank is also called the rank of $(M,\varphi)$. In this paper we will
use the shorter terminology (connected) Breuil module. 

\subsection{Classification results}\label{2.1}

We recall the following fundamental classification (see \cite{K1}, Theorem 2.3.5, \cite{VZ1}, Theorem 1, and \cite{Lau1}, Theorem 1.1 for $p>2$ and see \cite{K2}, \cite{Lau2}, Corollary 7.7, \cite{Liu3}, Theorem 1.0.2, and \cite{Kim}, Corollary 4.3 for $p=2$).

\begin{thm}\label{T3}
  There exists a contravariant functor $\mathbb{B}:\mathcal
  F\rightarrow\mathcal B$ which is an antiequivalence of categories, which is
  $\mathbb{Z}_p$-linear, and which takes short exact sequences (in the
  category of abelian sheaves in the faithfully flat topology of
  $\Spec V$) to short exact sequences (in the category of
  $\mathfrak{S}$-modules endowed with Frobenius maps).
\end{thm}

It is easy to see that $\mathbb{B}$ induces an antiequivalence of
categories $\mathbb{B}:\mathcal F_1\rightarrow\mathcal B_1$ which takes short
exact sequences to short exact sequences. For an object $G=\Spec R$ of $\mathcal F$, let $o(G)\in \mathbb{N}$ be such that $p^{o(G)}$ is the order of $G$ i.e., is the rank of $R$ over $V$.

We check that if $G$ is an object of $\mathcal F_1$, then the object
$\mathbb{B}(G)$ of $\mathcal B_1$ has rank $o(G)$. The contravariant
Dieudonn\'e module $\mathbb{D}(G_k)$ of $G_k$ is equal to
$k\otimes_{\sigma,\mathfrak{S}}\mathbb{B}(G)$, where we denote also by
$\sigma$ its composite with the epimorphism
$\mathfrak{S}\twoheadrightarrow k=\mathfrak{S}/(p,u)$ (the covariant
part of this statement follows from either \cite{Z1}, Theorem 6 or
\cite{Z2}, Theorem 1.6 provided $G$ is connected and from \cite{Lau2}, Section 8 in general, once we recall that $G$ is the kernel of an
isogeny of $p$-divisible groups over $V$). This implies that
the rank of $\mathbb{B}(G)$ is $o(G)$.

Let $H$ be an object of $\mathcal F$. If $p^n$ annihilates $H$, then
to the chain of natural epimorphisms
$$H\twoheadrightarrow H/H\lceil p\rceil\twoheadrightarrow H/\lceil p^2\rceil\twoheadrightarrow \cdots \twoheadrightarrow H/H\lceil p^n\rceil=0$$
corresponds a normal series of the Breuil module
$\mathbb{B}(H)=(M,\varphi)$
$$0=(M_n,\varphi_n)\subseteq (M_{n-1},\varphi_{n-1})\subseteq \cdots \subseteq (M_0,\varphi_0)=(M,\varphi)$$ 
by Breuil submodules whose quotient factors are objects of $\mathcal
B_1$.  As each $M_{i-1}/M_i$ is a free $\mathfrak{S}_1$-module of
finite rank, the multiplication by $u$ map $u:M\rightarrow M$ is
injective. One computes the order $p^{o(H)}$ of $H$ via the formulas
$$o(H)=o(M,\varphi):=\sum_{i=1}^n \text{rank}_{\mathfrak{S}_1}(M_{i-1}/M_i)=\text{length}_{\mathfrak{S}_{(p)}}(M_{(p)}).$$

\begin{prop}\label{P1} 
  Let $f:G\rightarrow H$ be a morphism of $\mathcal F$. We write
  $g:=\mathbb{B}(f):\mathbb{B}(H)=(M,\varphi)\rightarrow
  \mathbb{B}(G)=(N,\psi)$. Then we have:

  \medskip {\bf (a)} The homomorphism $f_K:G_K\rightarrow H_K$ is a closed
  embedding if and only if the cokernel of $g:M\rightarrow N$ is annihilated
  by some power of $u$.

  \smallskip {\bf (b)} The homomorphism $f_K:G_K\rightarrow H_K$ is an
  epimorphism if and only if the $\mathfrak{S}$-linear map $g:M\rightarrow N$
  is a monomorphism.

  \smallskip {\bf (c)} The homomorphism $f_K:G_K\rightarrow H_K$ is an
  isomorphism if and only if the $\mathfrak{S}$-linear map $g:M\rightarrow N$
  is injective and its cokernel is annihilated by some power of $u$.
\end{prop} {\bf Proof.} Let $\tilde N:=\Coker(g)$. We prove (a). We
first show that the assumption that $f_K$ is not a closed embedding
implies that $\tilde N$ is not annihilated by a power of $u$. This
assumption implies that there exists a non-trivial flat closed
subgroup scheme $G_0$ of $G$ which is contained in the kernel of $f$
and which is annihilated by $p$. Let $(N_0,\psi_0):=\mathbb{B}(G_0)$;
the $\mathfrak{S}_1$-module $N_0$ is free of positive rank. From the
fact that $G_0$ is contained in the kernel of $f$ and from Theorem
\ref{T3}, we get that we have an epimorphism $\tilde
N\twoheadrightarrow N_0$. Thus $\tilde N$ is not annihilated by a
power of $u$.

To end the proof of part (a) it suffices to show that the assumption
that $\tilde N$ is not annihilated by a power of $u$ implies that
$f_K$ is not a closed embedding. Our assumption implies that also
$N_1:=\tilde N/p\tilde N$ is not annihilated by a power of $u$. As
$\mathfrak{S}_1=k[[u]]$ is a principal ideal domain, we have a unique
short exact sequence of $\mathfrak{S}_1$-modules
$$0\rightarrow N_2\rightarrow N_1\rightarrow N_0\rightarrow 0,$$ 
where $N_2$ is the largest $\mathfrak{S}_1$-submodule of $N_1$
annihilated by a power of $u$ and where $N_0$ is a free
$\mathfrak{S}_1$-module of positive rank. The $\mathfrak{S}$-linear
map $N_1^{(\sigma)}\rightarrow N_1$ induced naturally by $\psi$ maps
$N_2^{(\sigma)}$ to $N_2$ and thus it induces via quotients a
$\mathfrak{S}$-linear map $\psi_0:N_0^{(\sigma)}\rightarrow N_0$. As $\psi_0$
is a quotient map of $\psi$, its cokernel is annihilated by
$E$. Therefore the pair $(N_0,\psi_0)$ is a Breuil module. From this
and Theorem \ref{T3} we get that there exists a non-trivial flat
closed subgroup scheme $G_0$ of $G$ which is contained in the kernel
of $f$ and for which we have $\mathbb{B}(G_0)=(N_0,\psi_0)$. Thus
$f_K$ is not a closed embedding. Therefore (a) holds.

Part (b) is proved similarly to (a). Part (c) follows from (a) and
(b).\endproof

\subsection{Basic lemmas}

Let $n\in\mathbb{N}^{\ast}$. Let $H$ be a
truncated Barsotti--Tate group of level $n$ over $V$. Let $(M,\varphi):=\mathbb{B}(H)$;
thus $M$ is a free $\mathfrak{S}_n$-module of finite rank. Let $d$ be
the dimension of $H$ and let $h$ be the height of $H$. There exists a
direct sum decomposition $M=T\oplus L$ into free
$\mathfrak{S}_n$-modules such that the image of $\varphi$ is $ET\oplus
L$ and $T$ has rank $d$. The existence of such a direct sum
decomposition follows (via reduction modulo $p^n$) from the existence of the normal decomposition of
the Breuil window relative to $q_{\pi}:\mathfrak{S}\twoheadrightarrow V$ of any $p$-divisible group $\mathcal H$ over $V$ such that $\mathcal H[p^n]=H$ (for the covariant context with $p>2$ see \cite{VZ1}, Section 2; the
case $p=2$ is the same). The existence of the direct sum decomposition
$M=T\oplus L$ implies the existence of two $\mathfrak{S}_n$-bases
$\{e_1,\ldots,e_h\}$ and $\{v_1,\ldots,v_h\}$ of $M$ such that for
$i\in\{d+1,\ldots,h\}$ we have $\varphi(1\otimes e_i)=v_i$ and for
$i\in\{1,\ldots,d\}$ the difference $\varphi(1\otimes e_i)-Ev_i$
belongs to the $\mathfrak{S}$-submodule of $M$ generated by the
elements $v_{d+1},\ldots,v_h$. Indeed, we consider the composite map
\begin{displaymath}
  M \rightarrow M^{(\sigma)} \rightarrow \varphi( M^{(\sigma)})\rightarrow  M/T=L,
\end{displaymath}
where the first map $m \mapsto 1\otimes m$ is $\sigma$-linear, where
the second map is $\varphi$, and where the third map is induced by the
natural projection of $M$ on $L$ along $T$.  We tensor this composite
map with the residue field $k$ of $\mathfrak{S}_n$
\begin{displaymath}
  k \otimes_{\mathfrak{S}_n} M \rightarrow k \otimes_{\mathfrak{S}_n} M^{(\sigma)} \rightarrow k \otimes_{\mathfrak{S}_n} \varphi(M^{(\sigma)}) \rightarrow k \otimes_{\mathfrak{S}_n} M/T=k \otimes_{\mathfrak{S}_n}L.
\end{displaymath}
As $k$ is a perfect field, all these three maps are
surjective. Therefore we find by the lemma of Nakayama a
$\mathfrak{S}_n$-basis $\{e_1, \ldots, e_h\}$ of $M$ such that the
images $\varphi(1 \otimes e_{d+1}), \ldots, \varphi(1 \otimes e_h)$
form a $\mathfrak{S}_n$-basis of $M/T$. We can take $L$ to be the
$\mathfrak{S}_n$-submodule of $M$ generated by $v_{d+1} := \varphi(1
\otimes e_{d+1}), \ldots, v_h := \varphi(1 \otimes e_h)$. We choose an
arbitrary $\mathfrak{S}_n$-basis $\{v_1, \ldots, v_d\}$ of $T$. As the
image of $\varphi$ is $ET \oplus L$ we obtain the desired
$\mathfrak{S}_n$-basis $\{v_1,\ldots,v_h\}$ by making a suitable
change of the $\mathfrak{S}_n$-basis $\{v_1, \ldots, v_d\}$ of $T$.

\begin{lemma}\label{Lau1}
  Let $t\in\mathbb{N}^{\ast}$.  Let $x\in \frac{1}{u^t}M$ be such that
  $\varphi(1\otimes x)\in \frac{1}{u^t}M$. We write $x=\sum_{i=1}^h
  \frac{\alpha_i}{u^t} e_i$, where $\alpha_i\in \mathfrak{S}_n$. Then
  for each $i\in\{1,\ldots,h\}$ we have $E\sigma(\alpha_i)\in
  u^{t(p-1)} \mathfrak{S}_n$ (equivalently,
  $\frac{\sigma(\alpha_i)}{u^{tp}} E \in \frac{1}{u^t}\mathfrak{S}_n$)
  and thus $\alpha_i\in (p^{n-1},u)\mathfrak{S}_n$.
\end{lemma} {\bf Proof.}
We compute
\begin{displaymath}
  \varphi(1 \otimes x) = \sum_{i=1}^{d} \frac{\sigma(\alpha_i)}{u^{tp}} \varphi (1\otimes e_i) +
  \sum_{i=d+1}^{h} \frac{\sigma(\alpha_i)}{u^{tp}} \varphi (1\otimes e_i)
\end{displaymath}
\begin{displaymath}
  = \sum_{i=1}^{d} \frac{\sigma(\alpha_i)}{u^{tp}} Ev_i +
  \sum_{i=d+1}^{h} \eta_i v_i \in \frac{1}{u^t}M,
\end{displaymath}
for suitable elements $\eta_i \in \frac{1}{u^{t}}\mathfrak{S}_n$.  For
$i\in\{1, \ldots, d\}$ this implies directly that
\begin{equation}\label{F1}
  \frac{\sigma(\alpha_i)}{u^{tp}} E \in \frac{1}{u^t}\mathfrak{S}_n.
\end{equation} 
The $\eta_i$ for $i\in\{d+1, \ldots, h\}$ are of the form
\begin{displaymath}
  \eta_i = \frac{\sigma(\alpha_i)}{u^{tp}} + 
  \sum_{j=1}^{d} \lambda_j \frac{\sigma(\alpha_j)}{u^{tp}},
\end{displaymath}
for some elements $\lambda_j \in \mathfrak{S}_n$. If we multiply the
last equation by $E$ we obtain from (\ref{F1}) that its belonging
relation also holds for $i\in\{d+1,\ldots, h\}$. For $i\in \{1,\ldots, h\}$, we get that $\sigma(\alpha_i)E(0)=\sigma(\alpha_i(0))a_0=0\in W_n(k)$ and this implies that $\alpha_i(0)\in p^{n-1}W_n(k)$; therefore we have $\alpha_i\in (p^{n-1},u)\mathfrak{S}_n$.\endproof

\begin{lemma}\label{Lau2}
  Let $t\in\mathbb{N}$. Let $N$ be a $\mathfrak{S}_n$-submodule of $\frac{1}{u^t}M$ which
  contains $M$. We assume that $\varphi$ induces a
  $\mathfrak{S}$-linear map $N^{(\sigma)}\rightarrow N$. Then we have
  $p^tN\subseteq M$.

\end{lemma} {\bf Proof.} We will prove this by induction on $t\in
\mathbb{N}$. The case $t=0$ is trivial. For the passage from $t-1$ to
$t$ we can assume that $t>0$. Let $x\in N$. From Lemma \ref{Lau1}
applied to $x$ we get that $px\in u\frac{1}{u^t}M=
\frac{1}{u^{t-1}}M$. This implies that $pN\subseteq
\frac{1}{u^{t-1}}M$. Let $\tilde N:=pN+M\subseteq
\frac{1}{u^{t-1}}M$. It is easy to see that $\varphi$ induces a
$\mathfrak{S}$-linear map $\tilde N^{(\sigma)}\rightarrow \tilde N$. By
induction applied to $\tilde N$ we get that $p^{t-1}\tilde N\subseteq
M$. This implies that $p^tN\subseteq M$. This ends the
induction.\endproof

\medskip The technical part of the method we use in this paper can be
summarized as follows. With $t$ and $N$ as in Lemma \ref{Lau2}, we will vary the uniformizer $\pi$ of $V$ to obtain a universal sharp
upper bound $t_0$ of $t$ (i.e., to show that in fact we have the first universal inclusion $N\subseteq \frac{1}{u^{t_0}}M$) and to refine Lemma \ref{Lau2} so that we get a smaller number $s\in\{0,\ldots,t_0\}$ for which we have the second universal inclusion $p^sN\subseteq M$.

 \section{Motivation, formulas, and bounds for $s$}
 In this section we present formulas for $s$ and upper
 bounds of it in terms of $e$. The main technical result that lies
 behind these formulas is also presented in this section (see
 Proposition \ref{P2} below). 

For a real number $x$, let $\lfloor x\rfloor$ be the integral part (floor) of $x$ (i.e., the greatest integer smaller or equal to $x$). We define
$$m:=\ord_p(e).$$ 
Let $a_e:=1$ and
\begin{displaymath}
  E_0 :=  \sum_{i\in p\mathbb{N}\cap [0,e]} a_iu^i= a_{p\lfloor \frac{e}{p}\rfloor}u^{p\lfloor \frac{e}{p}\rfloor}
  +\cdots +  a_{2p}u^{2p} + a_pu^p + a_0 \in W(k)[u^p].
\end{displaymath}
Let $E_1:=E-E_0\in W(k)[u]$. We define the numbers $\tau$ and $\iota$
as follows.

If $m=0$, then let $\tau(\pi):=1$ and $\iota(\pi):=0$.

If $m\ge 1$, let $\tau(\pi)\in\mathbb{N}^{\ast}\cup\{\infty\}$ be the
content $\ord_p(E_1)$ of $E_1$. Thus
$$\tau(\pi):=\ord_p(E_1)=\min\{\ord_p(a_i)|i\in\{1,2,\ldots,e-1\}\setminus 
p\mathbb{N}^{\ast}\}.$$ If $m\ge 1$ and $\tau(\pi)<\infty$, let
$\iota(\pi)\in \{1,2,\ldots,e-1\}\setminus p\mathbb{N}^{\ast}$ be the
smallest number such that we have
$$\tau(\pi)=\ord_p(a_{\iota(\pi)}).$$
\indent For all $m\ge 0$ we define
$$\tau=\tau_V:=\min\{\tau(\pi)|\pi\;\;\text{a}\;\;\text{uniformizer}\;\;\text{of}\;\;V\}.$$
If $\tau<\infty$, then we also define
$$\iota=\iota_V:=\min\{\iota(\pi)|\pi\;\;\text{a}\;\;\text{uniformizer}\;\;\text{of}\;\;V\;\;\text{with}\;\;\tau(\pi)=\tau\}.$$

\begin{lemma}\label{L3}
  We have $\tau\le m+1<\infty$.
\end{lemma} {\bf Proof.} If $m=0$, then this holds as $\tau=1$. Thus
we can assume that $m\ge 1$. We consider a new uniformizer
$\tilde\pi:=\pi+p$ of $V$. The unique Eisenstein polynomial $\tilde
E(u)=u^e+\tilde a_{e-1}u^{e-1}+\cdots+\tilde a_1u+\tilde a_0$ in $u$
with coefficients in $W(k)$ that has $\tilde \pi$ as a root is $\tilde
E(u)=E(u-p)$. Thus $\tilde
a_{e-1}=-pe+a_{e-1}=-p^{m+1}e^\prime+a_{e-1}$, where $e^\prime:=p^{-m}e\in\mathbb{N}^{\ast}\setminus p\mathbb{N}^{\ast}$. Therefore $p^{m+2}$
does not divide either $a_{e-1}$ or $\tilde a_{e-1}$. This implies
that either $\tau(\pi)\le m+1$ or $\tau(\tilde\pi)\le m+1$. Thus
$\tau\le m+1<\infty$.\endproof
  
\begin{prop}\label{P2}
  Let $n$ and $t$ be positive integers. Let $C=C(u)\in\mathfrak{S}$ be
  a power series whose constant term is not divisible by $p^n$. We
  assume that
  \begin{equation}\label{F2}
    E \sigma(C) \in (u^t,p^n)\mathfrak{S}.
  \end{equation}
If $\tau(\pi) = \infty$, then we have $t \le ne$. If $\tau (\pi) \neq \infty$, then we have
$$t \leq \min\{\tau(\pi)e+\iota(\pi),ne\}.$$
Moreover, if $m=0$, then we have $p\sigma(C) \in (u^t,p^n)\mathfrak{S}$
and if $m\geq 1$, then 
  we have $p^{\tau(\pi)+1}\sigma(C) \in (u^t,p^n)\mathfrak{S}$ (if $m\ge 1$ and the content of $C$ is $1$, then in fact we have $\tau(\pi)+1\ge n$).
\end{prop} {\bf Proof.} Clearly we can remove from $C$ all monomials
of some degree $i$ that satisfies the inequality $pi \geq t$. Therefore we can assume that $C$
is a polynomial of degree $d$ such that $pd < t$. 
As $E_0\sigma(C)$ and $E_1\sigma(C)$ do not contain monomials of the
same degree, the relation $E \sigma(C) \in (u^t,p^n)\mathfrak{S}$
implies that
\begin{displaymath}
  E_0 \sigma(C) \in (u^t,p^n)\mathfrak{S} \quad \text{and} \quad 
  E_1 \sigma(C) \in (u^t,p^n)\mathfrak{S}.
\end{displaymath}

We first consider the case when $p$ does not divide $e$ (i.e.,
$m=0$). Thus $E_0$ divided by $p$ is a unit of
$\mathfrak{S}$. Therefore we have $p \sigma(C) \in
(u^t,p^n)\mathfrak{S}$. As $p\deg(C)=pd<t$, we get that $p \sigma(C)
\equiv 0 \, \mod p^n$. As moreover $E_1 \equiv u^e \, \mod p$, we get $u^e
\sigma(C) \in (u^t,p^n)\mathfrak{S}$.  As the constant term of $C$ is
not divisible by $p^n$ this implies that $t \leq
e=\min\{\tau(\pi)e+\iota(\pi),ne\}$.

From now on we will assume that $p|e$. By the
Weierstra{\ss} preparation theorem (\cite{Bou}, Chapter 7, Section 3,
number 8) we can assume that $C$ is a Weierstra{\ss} polynomial of
degree $d$ (i.e., a monic polynomial of degree $d$ such that $C-u^d$
is divisible by $p$) and thus of content $1$. Indeed, to check this let $c\in\{0,\ldots,n-1\}$ be such that $p^c$ is the content of $C$. We set $\bar{C} := (1/p^c)C$.  The constant term of $\bar{C}$ is not divisible by $p^{n-c}$. As $E \sigma(\bar{C}) \in
(u^t,p^{n-c})\mathfrak{S}$, it suffices to prove the proposition for
$(\bar{C},n-c)$ instead of $(C,n)$. But $\bar{C}$ is a unit times a Weierstra{\ss} polynomial.

Before we continue, we first prove the following lemma.

\begin{lemma}\label{L4}
  Let $n$ and $t$ be positive integers. We assume that $p$ divides
  $e$. Let
  \begin{displaymath}
    C =C(u)= u^d + c_{d-1}u^{d-1} + \cdots +c_1u + c_0 \in W(k)[u]
  \end{displaymath}
  be a Weierstra{\ss} polynomial such that $pd < t$ and $c_0 \notin
  p^nW(k)$.  We also assume that
  \begin{displaymath}
    E_0 \sigma (C) \in (u^t,p^n)\mathfrak{S}.
  \end{displaymath}

  Then $d = (n-1)\frac{e}{p}$ and for each $i\in\{0,1,\ldots,n-1\}$ we
  have:

\begin{equation}\label{F3}
  \ord_p(c_{i\frac{e}{p}}) = n-i-1, \quad \text{and} \quad 
  \ord_p(c_j) \geq n-i, \quad \text{for} \quad  0\le j < i\frac{e}{p}.
\end{equation}
Moreover we also have $t \leq ne$.
\end{lemma} {\bf Proof:}
For $l\in\{0,\ldots,d-1\}$, let $\gamma_l=\sigma(c_l)\in W(k)$. Thus 
\begin{displaymath}
  \sigma(C) = u^{dp} + \gamma_{d-1}u^{(d-1)p} + \cdots +\gamma_1u^{p} + \gamma_0.
\end{displaymath}
We define $c_d$ and $\gamma_d$ to be $1$. For $l<0$ or $l>d$, we
define $c_l=\gamma_l=0$. We have $\ord_p(\gamma_l) = \ord_p(c_l)$ for
all $l\in\mathbb{Z}$. Moreover we set
\begin{displaymath}
  E_0 \sigma (C) = \sum_{j=0}^{d+\frac{e}{p}} \beta_{jp}u^{jp}, \quad\quad \beta_{jp} \in W(k).
\end{displaymath}
By our assumption $\beta_{jp}$ is divisible by $p^n$ if $jp < t$ and
in particular if $j \leq d$. For $j\in\{0,\ldots,d+\frac{e}{p}\}$ we
have the identity
\begin{equation}\label{F4}
  \beta_{jp} = a_0\gamma_j + a_p\gamma_{j-1} + \cdots + a_e\gamma_{j-\frac{e}{p}}.
\end{equation}

By induction on $j\in\{0,\ldots,d\}$ we show that (\ref{F3})
holds. This includes the equality part of (\ref{F3}) if $j =
i\frac{e}{p}$. Our induction does not require that $i < n$. But of
course, the assumption that $i \geq n$ in the equality part of
(\ref{F3}) (resp. the assumptions that $j\leq d$ and $i \geq n+1$ in
the inequality part of (\ref{F3})), leads (resp. lead) to a
contradiction as the order of any $c_j$ can not be negative.

The case $j=0$ follows by looking at the constant term of
$E_0\sigma(C)$.  The passage from $j-1$ to $j$ goes as follows.  Let
us first assume that $(i-1)\frac{e}{p} < j < i\frac{e}{p}$ for some
integer $i \geq 1$ (we do not require $i\frac{e}{p} \leq d$).

We show that the assumption $\ord_p(c_j)=\ord_p(\gamma_j) < n-i$ leads to a
contradiction. Indeed, in this case the first term of the right hand side of (\ref{F4}) would
have order $\leq n-i$ but all the other terms of the right hand side of (\ref{F4}) would have order
strictly bigger than $n-i$.  Note for example that $j - \frac{e}{p} <
(i-1)\frac{e}{p}$ and therefore $\ord_p(\gamma_{j-\frac{e}{p}}) \geq n
- i + 1$ by induction assumption.  As $\ord_p(\beta_{jp})\geq n$ we
obtain a contradiction to (\ref{F4}).

Finally we consider the passage from $j-1$ to $j$ in the case when $j
= i\frac{e}{p}$. Then we use the equation
   
\begin{equation}\label{F5}
  \beta_{ie} = a_0\gamma_{i\frac{e}{p}} + a_p\gamma_{i\frac{e}{p}-1} + \cdots + 
  a_e\gamma_{(i-1)\frac{e}{p}}.
\end{equation}
By the induction assumption we have
$\ord_p(a_e\gamma_{(i-1)\frac{e}{p}})=n-i$. From the inequality part
of (\ref{F3}) we get the inequalities
\begin{displaymath}
  \ord_p(a_p\gamma_{i\frac{e}{p}-1})\ge n-i+1,\ldots,
  \ord_p(a_{e-p}\gamma_{1+(i-1)\frac{e}{p}})\ge n-i+1.
\end{displaymath}
As $\beta_{ie}=\beta_{jp}\in p^nW(k)$, it follows from (\ref{F5}) that
\begin{displaymath}
  \ord_p(a_0\gamma_{i\frac{e}{p}})= \ord_p(a_e\gamma_{(i-1)\frac{e}{p}})=n-i.
\end{displaymath}
Thus $\ord_p(\gamma_{i\frac{e}{p}})=n-i-\ord_p(a_0)=n-i-1$. This ends
our induction.

We check that $d = (n-1)\frac{e}{p}$. If $d$ is of the form
$i\frac{e}{p}$, then $0 = \ord_p(c_d) = n - i -1$ gives $i = n-1$ and
thus $d$ is as required. 
We are left to show that the assumption that $d$ is not
of this form, leads to a contradiction. Let $i\in\mathbb{N}^{\ast}$ be such that $(i-1)\frac{e}{p}<d < i\frac{e}{p}$. Then the inequality $ 0 =
\ord_p(c_d)\geq n - i$ gives $i \geq n$. But then
$\ord_p(c_{(i-1)\frac{e}{p}}) = n - i\leq 0$ contradicts the assumption that
$C$ is a Weierstra{\ss} polynomial and therefore that $c_{(i-1)\frac{e}{p}}
\in pW(k)$. 

As $p$ divides $e$ and $d = (n-1)\frac{e}{p}$, $E_0\sigma(C)\in (u^t,p^n)\mathfrak{S}$ is a monic polynomial of degree $e+pd=ne$ and thus we have $t\le ne$.\endproof

\begin{cor}\label{C5}
  With the notations of Lemma \ref{L4}, let
  $l\in\{0,1,\ldots,e-1\}$. Let $E_2 = E_2(u)=u^l + b_{l-1}u^{l-1} +
  \cdots +b_1u + b_0 \in W(k)[u]$ be a Weierstra{\ss} polynomial of
  degree $l$. If we have $E_2 \sigma(C) \in (u^t, p^n)\mathfrak{S}$,
  then $l\geq t$.
\end{cor} {\bf Proof:} If $n=1$, then $d=0$ and $C=c_0$ is a unit of
$W(k)$; thus $E_2\in (u^t, p^n) \mathfrak{S}$ and therefore $l\geq
t$. Thus we can assume that $n\ge 2$. We write:
\begin{displaymath}
  E_2\sigma(C) = \sum_{i=0}^{l+pd} \delta_iu^i.
\end{displaymath}
Let $q := \lfloor\frac{l}{p}\rfloor < \frac{e}{p}$. We have the equation:
\begin{equation}\label{F6}
  \delta_l = \gamma_0 + b_{l-p}\gamma_1 + \cdots + b_{l - qp}\gamma_q.
\end{equation}
For $i\in\{1,\ldots,q\}$ we have $\ord_p(\gamma_i)=\ord_p(c_i)\geq n-1$ (cf. Lemma
\ref{L4} and $n\ge 2$); therefore $\ord_p(b_{l-ip}\gamma_i)\geq n$. As $\ord_p
(\gamma_0)=\ord_p(c_0)= n-1$, from (\ref{F6}) we get $\ord_p (\delta_l)=n-1$. From
this and the assumption $E_2 \sigma(C) \in (u^t, p^n)\mathfrak{S}$, we
get $l\geq t$. 

We note that the inequality $l\geq t$ implies that actually the case $n\ge 2$ does not occur (cf. the sequence of relations $e>l\geq t>pd=(n-1)e$).
\endproof

\medskip We are now ready to prove Proposition \ref{P2} in the 
case when $p|e$ (i.e., $m\ge 1$) and the content of $C$ is $1$. It suffices to show that $\tau(\pi)+1\ge n$. We can assume that $\tau(\pi)+1\leq n$ and it suffices to show that $\tau(\pi)+1=n$.

As $\tau(\pi)=\ord_p(E_1)$, by Weierstra{\ss} preparation theorem we
can write
\begin{displaymath}
  E_1 = p^{\tau(\pi)}E_2\theta,
\end{displaymath}
where $\theta \in\mathfrak{S}$ is a unit and where $E_2\in W(k)[u]$ is
a Weierstra{\ss} polynomial of degree $\iota(\pi)< e$. As $E_1\sigma
(C) \in (u^t,p^n)\mathfrak{S}$ and as $n>\tau(\pi)$, we get that
\begin{displaymath}
  E_2\sigma (C) \in (u^t,p^{n-\tau(\pi)})\mathfrak{S}.
\end{displaymath}
As $d = (n-1)\frac{e}{p}\geq \tau(\pi)\frac{e}{p}$ and $c_d=1$, we can consider the monic polynomial
\begin{displaymath}
  C_1 = C_1(u)= u^d + c_{d-1}u^{d-1} + \cdots + c_{\tau(\pi)\frac{e}{p}}u^{\tau(\pi)\frac{e}{p}}\in W(k)[u]. 
\end{displaymath}
It follows from Lemma \ref{L4} that $\ord_p(c_j) \geq n - \tau(\pi)$
for $j < \tau(\pi)\frac{e}{p}$. Thus $C_1-C\in
p^{n-\tau(\pi)}\mathfrak{S}$ and therefore we obtain
\begin{equation}\label{F7}
  E_2\sigma (C_1)=E_2\sigma(C_1-C)+E_2\sigma(C)\in (u^t,p^{n-\tau(\pi)}) \mathfrak{S}.
\end{equation}
We write $C_1 = u^{\tau(\pi)\frac{e}{p}}C_2$. Then the constant term
of $C_2$ is $c_{\tau(\pi)\frac{e}{p}}$ and thus it is not divisible by
$p^{n-\tau(\pi)}$, cf. (\ref{F3}).  As $n>\tau(\pi)$, the
relation $(n-1)e = pd<t$ implies that $t- \tau(\pi)
e>0$. Thus from (\ref{F7}) we get that
\begin{equation}\label{F8}
  E_2 \sigma (C_2) \in (u^{t - \tau(\pi) e},p^{n-\tau(\pi)}) \mathfrak{S}.
\end{equation}
A similar argument shows that
\begin{equation}\label{F9}
  E_0 \sigma (C_2) \in (u^{t - \tau(\pi) e}, p^{n-\tau(\pi)}) \mathfrak{S}.
\end{equation}
From Corollary \ref{C5} applied to the quintuple $(t- \tau(\pi)
e,C_2,E_0,E_2,n - \tau(\pi))$ instead of $(t,C,E_0,E_2,n)$, we get that
$\iota(\pi)=\deg(E_2)\ge t- \tau(\pi) e$. As $\iota(\pi)\le e-1$ and as
$n\ge \tau(\pi)+1$, we conclude that $t\le \tau(\pi)
e+\iota(\pi)=\min\{\tau(\pi) e+\iota(\pi),ne\}$. As $\iota(\pi)\le e-1$, the
relations $(n-1)e = pd<t\le \min\{\tau(\pi) e+\iota(\pi),ne\}$ imply that
$n\le \tau(\pi)+1$. Thus $n=\tau(\pi)+1$ and this ends the proof of the Proposition \ref{P2} in the
second case $p|e$.\endproof

\subsection{Formulas for $s$}\label{REC} 
For a uniformizer $\pi$ of $V$, Lemma \ref{Lau1} and Proposition \ref{2}
motivate the introduction of the following invariant
\begin{displaymath}
  t(\pi) := \lfloor \frac{\tau (\pi)e + \iota(\pi)}{p-1} \rfloor\in\mathbb{N}\cup\{\infty\}.
\end{displaymath}
Always there exists a $\pi$ such that $t(\pi)$ is finite, cf. Lemma
\ref{L3}.

Based on the last sentence of Proposition \ref{2}, we define
$\epsilon\in\{0,1\}$ as follows. If $m=0$, then $\epsilon:=0$ and if
$m\ge 1$, then $\epsilon:=1$. Next we will define a number
$s=s_V\in\mathbb{N}$ as follows. For all non-negative integers $i$ let
$$s_i:=i(\tau+\epsilon)\;\;\;\text{and}\;\;\;t_i:=\lfloor \frac{\tau e + \iota}{(p-1)p^i} \rfloor.$$
Thus $t_0=\lfloor \frac{\tau e + \iota}{p-1} \rfloor$ is the minimum of $t(\pi)$ for all possible uniformizers $\pi$ of $V$ and we have $t_{i+1}=\lfloor \frac{t_i}{p} \rfloor$ for all $i\geq 0$. Let $z\in\mathbb{N}$ be the smallest number such that we have $t_z-t_{z+1}\le\tau+\epsilon$. Equivalently, let $z\in\mathbb{N}$ be the smallest number such that we have 
$$s_z+t_z=\min\{s_i+t_i|i\in\mathbb{N}\}=:s=s_V.$$

We note that
$$0=s_0<s_1<\cdots<s_z\;\;\;\text{and}\;\;\;0\le t_z<t_{z-1}<\cdots<t_0.$$ 
Moreover the following relations hold:
\begin{equation}\label{F10}
  s=s_z+t_z< s_{z-1}+t_{z-1}< \cdots < s_1+t_1<s_0+t_0=t_0=\lfloor \frac{\tau e + \iota}{p-1} \rfloor.
\end{equation}

By Lemma \ref{L3} there exists a $\pi$ such that $\tau(\pi)\le m+1$
(with $m = \ord_p(e)$).  We have $s\le \lfloor \frac{\tau e
  + \iota}{p-1} \rfloor$ and thus, as $\iota \leq e-1$, we get
\begin{equation}\label{F11}
  s\le \frac{2e-1+e\ord_p(e)}{p-1}.
\end{equation} 

If $t_0\geq 1$, let $v:=\lfloor \log_p t_0\rfloor$ and we consider the $p$-adic expansion
\begin{equation}\label{F13exp}
  t_0 = \sum_{j=0}^{v} a_jp^j, \quad 0\leq a_j < p, \; a_v \neq 0.
\end{equation}

\noindent
{\bf Example 1}
We assume that $e\le p-2$. Thus $p$ is odd, $m=\iota=0$, and
$\tau=1$. We have $s\le \frac{e}{p-1}$, cf. (\ref{F10}).  Therefore
$s=0$.

\medskip\noindent {\bf Example 2} We assume that $m=0$ and $e\geq p-1$. As $m=0$ we have $\tau = \tau + \epsilon = 1$ and thus $(s_z,t_z)$ is either $(v,1)$ or $(v+1,0)$. Therefore for $m=0$ and $e\geq p-1$ we get $s = v
+1=1+\lfloor \log_p({{e}\over {p-1}})\rfloor=\lfloor \log_p({{pe}\over {p-1}})\rfloor$ and this is the same expression as in \cite{Bon}. We have $s = 1$ if and only if $p-1 \leq e \leq p^2 - p - 1$. In general, for $m=0$ we have the following practical upper bound $s \leq 1+\log_p e$.

\medskip\noindent {\bf Example 3} We assume that $E=u^p-p$. Then
$\tau(\pi)=\infty$. For $n\in \mathbb{N}^{\ast}$ we have
$E\sigma(\sum_{i=1}^n p^{n-i}u^{i-1})=u^{pn}-p^n$. Thus if $M=\mathfrak{S}_n$ and $\varphi:M^{(\sigma)}\to M$ takes $1$ to the image of $E$ in $M=\mathfrak{S}_n$, then for $N:=M+{1\over {u^n}}(\sum_{i=1}^n p^{n-i}u^{i-1})M$ we have $\varphi(N^{(\sigma)})=M\subset N$, $p^nN=0\subset M$, $p^{n-1}N\nsubseteq M$ and therefore Lemma \ref{Lau2} is optimal in general. The smallest $t$ for which one has $N\subset {1\over {u^t}}M$ is $n$ itself and thus it does not have an upper bound independently of $n$.

Let $\tilde \pi=c_0p+c_1\pi+\cdots +c_{p-1}\pi^{p-1}$ be another
uniformizer of $V$. We have $c_0,c_1,c_1^{-1},c_2,\ldots,c_{p-1}\in
W(k)$. It is easy to see that $\tilde\pi^p$ is congruent modulo $p^2$
to the element $pc_1^p$ of $pW(k)$. Thus, if $\tilde E(u)=u^p+\tilde
a_{p-1}u^{p-1}+\cdots+\tilde a_0$ is the Eisenstein polynomial that
has $\tilde \pi$ as a root, then for $i\in\{1,\ldots,p-1\}$ we have
$\tilde a_i\in p^2W(k)$. Therefore $\tau\geq 2$. But $\tau\le m+1=2$,
cf. Lemma \ref{L3}. Thus $\tau=2$. From inequalities (\ref{F10}) we get
that we have $s\le \frac{2p+\iota}{p-1}$. Thus for $p=3$ we
have $s\le 4$ and for $p\ge 5$ we have $s\le 3$.

\medskip\noindent {\bf Example 4} We assume that $m \geq 1$. We have $t_{v+1}\leq {{t_0}\over {p^{v+1}}}<1$ and thus $t_{v+1}=0$. Therefore $z\le v+1$ and $s \leq s_{v+1}+t_{v+1}=s_{v+1}=(v+1)(\tau + 1)$. We have the inequalities $(\tau + 1)e/(p-1) \geq t_0$ and $\log_p t_0 \geq v$.  Moreover $\tau \leq m+1$ by Lemma \ref{L3}. We find the inequality
\begin{displaymath}
   s \leq (\log_p (m+2) + \log_p e - \log_p(p-1) + 1)(m+2).
\end{displaymath}
As $m = \ord_p e$, $\log_p(m+2) < m+1$, and $\log_p(p-1)\in [0,1)$, we find that
\begin{displaymath}
   s< (\log_p e + \ord_p e +2)(\ord_p e +2).
\end{displaymath}
If $(p,m)\neq (2,1)$, then $\log_p(m+2) \leq m$ and we get that
\begin{displaymath}
   s\leq (\log_p e + \ord_p e +1)(\ord_p e +2).
\end{displaymath}

\section{Proof of Theorem \ref{T1}}

We prove that Theorem \ref{T1} holds for the number $s$ of Subsection
\ref{REC} if $V$ is complete and has perfect residue field $k$. We
already noted in Section 2 that it is enough to treat this case. 

We choose the uniformizer $\pi$ of $V$ in such a way that $\tau =
\tau(\pi)$ is minimal and $\iota=\iota(\pi)$. Let $z\in\mathbb{N}$, the
sequence of pairs $(s_0,t_0),\ldots,(s_z,t_z)$, and $s\in\mathbb{N}$ be as
in Subsection \ref{REC}. Thus $ t_0 = \lfloor (\tau e + \iota)/(p-1)
\rfloor=t_0(\pi)$ and $s=s_z+t_z$.

Let $f: G \rightarrow H$ be a homomorphism of finite flat commutative
group schemes of $p$ power order over $V$, which induces an
isomorphism $f_K:G_K\rightarrow H_K$ in generic fibers. We first remark that
if there exists a homomorphism $f':H\rightarrow G$ such that $f\circ
f'=p^s\id_H$, then we have $f'_K=p^sf_K^{-1}$ and thus the
equality $f'\circ f=p^s\id_G$ also holds as it holds
generically; moreover, as the equalities continue to hold in the special
fibre we conclude that the kernel and cokernel of $f_k$ are annihilated
by $p^s$.

We choose an epimorphism $\xi_H:\tilde{H} \twoheadrightarrow H$ from a
truncated Barsotti--Tate group $\tilde{H}$. We consider
the following fiber product
\begin{displaymath}
  \begin{CD}
    \tilde{G} @>{\tilde{f}}>> \tilde{H}\\
    ^{\xi_G}@VVV   @VVV^{\xi_H}\\
    G @>{f}>> H
  \end{CD}
\end{displaymath}
in the category $\mathcal{F}$. Then $\tilde{f}_K$ is an isomorphism. 
Assume that there exists a homomorphism $\tilde{f}': \tilde{H} \rightarrow \tilde{G}$ 
such that $\tilde{f} \circ \tilde{f}' = p^s\id_{\tilde H}$. Then $\xi_G
\circ \tilde{f}'$ is zero on the finite flat group scheme $\Ker(\xi_H)$
because this is true for the generic fibers. Thus there exists $f': H
\rightarrow  G$ such that $f'\circ \xi_H = \xi_G \circ
\tilde{f}'$. One easily verifies that $f\circ f' = p^s\id_H$.

Thus to prove the existence of $f'$ we can assume that $f=\tilde{f}$
and that $H=\tilde{H}$ is a truncated Barsotti--Tate group of level
$n>s$.

Next we translate the existence of $f'$ in terms of Breuil
modules. Let $(M,\varphi)$ and $(N,\psi)$ be the Breuil modules of $H$
and $G$ (respectively). We know by Proposition \ref{P1} (c) that to
$f$ corresponds a $\mathfrak{S}$-linear monomorphism $M
\hookrightarrow N$ whose cokernel is annihilated by some power
$u^t$. We will assume that $t$ is the smallest natural number with
this property. We put aside the case $t = 0$ (i.e., the case when
$f:G\rightarrow H$ is an isomorphism) which is trivial. The existence of
$f':H\rightarrow G$ is equivalent to the inclusion 

\begin{equation}\label{F13}
p^sN \subseteq M.   
\end{equation}

Before we prove this inclusion we show that Theorem \ref{T1} follows
from it.  It remains to prove the last sentence of the Theorem \ref{T1}. Hence
again let $H$ be a truncated Barsotti--Tate group of level $n > s$.

The identity $f\circ f'=p^s\id_H$ means that we have a
commutative diagram: 

\begin{displaymath}
  \xymatrix{H \ar[r]^{f'} \ar[dr]_{p^s} & G \ar[d]^f\\
    &  H.
  }
\end{displaymath}
We note that a homomorphism of finite flat group schemes over $V$ is zero if
it induces the zero homomorphism at the level of generic fibers. Therefore we obtain
a commutative diagram:
\begin{displaymath}
  \xymatrix{H/(H[p^s]) \ar[r]^{\quad \breve{f}'} \ar[dr] & G \ar[d]^f\\
    &  H.
  }
\end{displaymath}
We see that $\breve{f'}$ is a closed immersion because the diagonal arrow is. 
Now we apply the functor $\lceil p^{n-s}\rceil$ (see Section 1) 
to the last diagram:
\begin{displaymath}
  \xymatrix{H/(H[p^s])\; \ar[r]^{\; \breve{f}'\lceil p^{n-s} \rceil} \ar[dr] & G\lceil p^{n-s} \rceil 
    \ar[d]^{f\lceil p^{n-s}\rceil}\\
    &  H[p^{n-s}].
  }
\end{displaymath}
The horizontal homomorphism is again a closed immersion. As it is a homomorphism
between finite flat group schemes of the same order, it has to be an
isomorphism. The diagonal arrow is trivially an isomorphism. Thus
$f\lceil p^{n-s}\rceil: G\lceil p^{n-s} \rceil\rightarrow H[p^{n-s}]$ is an
isomorphism and therefore $G\lceil p^{n-s} \rceil$ is a truncated Barsotti--Tate
group of level $n-s$.  This shows the last sentence of Theorem \ref{T1}. It
remains to prove the inclusion (\ref{F13}).

We will prove by induction on $j\in\{0,\ldots,z\}$ that we have
$p^{s_j}N\subseteq \frac{1}{u^{t_j}}M$. 

We remark that by Lemma
\ref{Lau2} this implies that $p^{s_j + t_j}N \subset M$. As $s = s_z+ t_z$ the induction gives the desired inclusion (\ref{F13}) and ends the proof of Theorem \ref{T1}.
We also remark that already the base of the induction $j=0$
implies the Theorem \ref{T1} but with the much weaker bound $s = t_0$.

For the base of the induction it suffices to show that $t\le t_0$. Let
$x \in N$ be such that $u^{t-1}x \notin M$. With $\{e_1,\ldots,e_h\}$
a $\mathfrak{S}_n$-basis of $M$ as before Lemma \ref{Lau1}, we write
\begin{displaymath}
  x = \sum_{i=1}^h \frac{\alpha_i}{u^t} e_i, \quad 
  \alpha_i \in \mathfrak{S}_n.
\end {displaymath}  
From Lemma \ref{Lau1} we get that for each $i\in\{1,\ldots,h\}$ we have $E\sigma(\alpha_i) \in 
u^{t(p-1)}\mathfrak{S}_n$. By the minimality of $t$, there exists $i_0\in\{1,\ldots,h\}$ such that $\alpha_{i_0}$ is not divisible by $u$. Let $C=C(u)\in \mathfrak{S}$ be such that its reduction modulo $p^n$ is $\alpha_{i_0}$. The constant term of $C$ is not divisible by $p^n$ and we have $E\sigma(C) \in 
(u^{t(p-1)},p^n)\mathfrak{S}$. From this and Proposition \ref{P2}
we get that $t(p-1)\leq \min\{\tau e + \iota,ne\}$. Thus $t\le t_0=\lfloor (\tau e + \iota)/(p-1) \rfloor$. 

If $0<j<z$, then the inductive step from $j-1$ to $j$ goes as follows. We assume that  $p^{s_{j-1}}N\subseteq \frac{1}{u^{t_{j-1}}}M$. Let $l_{j-1}\in\{0,\ldots,t_{j-1}\}$ be the smallest number such that we have $p^{s_{j-1}}N\subseteq \frac{1}{u^{l_{j-1}}}M$. If $l_{j-1}=0$, then $p^{s_{j-1}}N\subseteq M$ and thus, as $s_{j-1}<s_j$, we also have $p^{s_j}N\subseteq M\subseteq \frac{1}{u^{t_j}}M$. Therefore we can assume that $1\le l_{j-1}\leq t_{j-1}$. Let $y\in p^{s_{j-1}}N$.  We write 
\begin{displaymath}
  y = \sum_{i=1}^h \frac{\eta_i}{u^{n_i}} e_i,
\end {displaymath}  
where $\eta_i \in \mathfrak{S}_n\setminus u\mathfrak{S}_n$ and where
$n_i\in\{0,\ldots,l_{j-1}\}$. Let $C_i=C_i(u)\in \mathfrak{S}$ be such
that its reduction modulo $p^n$ is $\eta_i$. We want to show that
$p^{\tau+\epsilon}y\in \frac{1}{u^{t_j}}M$. For this it suffices to
show that for each $i\in\{1,\ldots,h\}$ we have
$p^{\tau+\epsilon}\frac{\eta_i}{u^{n_i}}\in
\frac{1}{u^{t_j}}\mathfrak{S}_n$.  To check this we can assume that
$n_i\geq t_j+1$. As  

\begin{displaymath}
t_j+1=\lfloor \frac{t_{j-1}}{p}\rfloor+1 \geq \frac{t_{j-1}+1}{p}\geq
\frac{l_{j-1}+1}{p},
\end{displaymath}
we get that $pn_i-l_{j-1}\geq 1$. As $\frac{\eta_i}{u^{n_i}}=\frac{\eta_iu^{l_{j-1}-n_i}}{u^{l_{j-1}}}$, from Lemma \ref{Lau1} (applied with $y,\varphi(1\otimes y)\in p^{s_{j-1}}N\subseteq {{1}\over {u^{l_{j-1}}}}M$) we get that $E\sigma(C_iu^{l_{j-1}-n_i})\in (u^{(p-1)l_{j-1}},p^n)\mathfrak{S}$. This implies that 
$$E\sigma(C_i)\in (u^{pn_i-l_{j-1}},p^n)\mathfrak{S}\subseteq (u,p^n)\mathfrak{S}.$$ 
As $\eta_i \in \mathfrak{S}_n\setminus u\mathfrak{S}_n$, the constant term of $C_i$ is not divisible by $p^n$. Thus from Proposition \ref{P2} applied to $(C_i,pn_i-l_{j-1})$ instead of $(C,t)$ and from the definition of $\epsilon$ in Subsection \ref{REC}, we get $\sigma(p^{\tau+\epsilon}C_i)=p^{\tau+\epsilon}\sigma(C_i)\in (u^{pn_i-l_{j-1}},p^n)\mathfrak{S}$. This implies that we can write $p^{\tau+\epsilon}C_i=A_i+B_i$, where $A_i\in p^n\mathfrak{S}$ and where $B_i\in u^{n_i-\lfloor\frac{l_{j-1}}{p}\rfloor}\mathfrak{S}$. 
Thus
$$p^{\tau+\epsilon}\frac{\eta_i}{u^{n_i}}\in \frac{1}{u^{\lfloor\frac{l_{j-1}}{p}\rfloor}}\mathfrak{S}_n\subseteq \frac{1}{u^{\lfloor\frac{t_{j-1}}{p}\rfloor}}\mathfrak{S}_n=\frac{1}{u^{t_j}}\mathfrak{S}_n.$$ Therefore $p^{\tau+\epsilon}y\in \frac{1}{u^{t_j}}M$. As $s_j=s_{j-1}+\tau+\epsilon$ and as $y\in p^{s_{j-1}}N$ is arbitrary, we conclude that $p^{s_j}N=p^{\tau+\epsilon}p^{s_{j-1}}N\subseteq  \frac{1}{u^{t_j}}M$. This ends the induction.\endproof

\section{Proofs of Corollaries \ref{C1} to \ref{C4} and Theorem \ref{T2}}

\subsection{Proofs of the Corollaries  \ref{C1} to \ref{C4}}
Corollary \ref{C1} follows from  Theorem \ref{T1} and Example 1.

In connection to the Corollaries \ref{C2} to \ref{C4}, let
$\tilde{G}$ be the schematic closure in $G\times_V H$ of the graph of
the homomorphism $h:G_K\rightarrow H_K$ and let $i:\tilde{G}\hookrightarrow G\times_V H$ be the resulting closed embedding homomorphism. Via the two projections $q_1:G\times_V
H\rightarrow G$ and $q_2:G\times_V H\rightarrow H$ we get homomorphisms $\rho_1:\tilde
G\rightarrow G$ and $\rho_2:\tilde G\rightarrow H$. The generic fiber $\rho_{1,K}$ of
$\rho_1$ is an isomorphism.  

We prove the first part of Corollary \ref{C2}. We consider the commutative diagram:

\begin{displaymath}
   \xymatrix{
\tilde{G} \ar[r]^-i \ar[rd]_{\rho_1} & G \times_V H \ar[r]^-{q_2}
\ar[d]^{q_1} & H\\ 
& G. &
}
\end{displaymath}
By Theorem \ref{T1} there exists a homomorphism $\rho_1': G \rightarrow \tilde{G}$ such
that $\rho_1' \circ \rho_1 = p^s\id_{\tilde{G}}$. Then $q_2 \circ i
\circ \rho_1'$ is the desired extension of $p^s h.$ 

To check the last part of Corollary \ref{C2}, let $0\to G\to J\to H\to 0$ be a short exact sequence whose generic fibre splits. It defines an arbitrary element $\nu\in \Ker(\Ext^1(H,G)\to \Ext^1(H_K,G_K))$. Let $h:H_K\to J_K$ be a homomorphism that is a splitting of $0\to G_K\to J_K\to H_K\to 0$. Let $g:H\to J$ be such that its generic fibre is $p^sh$, cf. first part of Corollary \ref{C2}. Let $0\to G\to J_s\to H\to 0$ be the pull back of $0\to G\to J\to H\to 0$ via $p^s\id_H$. Then there exists a unique section $g_s:H\to J_s$ of $0\to G\to J_s\to H\to 0$ whose composite with $J_s\to J$ is $g$. Thus $p^s\nu=0$. This proves Corollary \ref{C2}.

We prove Corollary \ref{C3}. From Corollary \ref{C2} we get that
$p^sh$ extends to a homomorphism $G\rightarrow H$ which induces a homomorphism $G[p^{n-s}]=G/G[p^s]\rightarrow H[p^{n-s}]$ whose generic fibre is $h[p^{n-s}]$. This proves Corollary \ref{C3}.  

We prove Corollary \ref{C4} for $n>2s$. As $h:G_K\rightarrow H_K$ is an isomorphism, $\rho_{2,K}$ is also an isomorphism. Thus $\rho_2\lceil p^{n-s}\rceil:\tilde G\lceil p^{n-s}\rceil\rightarrow H[p^{n-s}]$ is an isomorphism, cf. Theorem \ref{T1}. By applying Theorem \ref{T1} to the Cartier dual $(\rho_1\lceil p^{n-s}\rceil)^{\text{t}}:( G\lceil p^{n-s}\rceil)^{\text{t}}\rightarrow (\tilde G\lceil p^{n-s}\rceil)^{\text{t}}$ of $\rho_1\lceil p^{n-s}\rceil:\tilde G\lceil p^{n-s}\rceil\rightarrow G\lceil p^{n-s}\rceil$, we get that $( G\lceil p^{n-s}\rceil)^{\text{t}}\lceil p^{n-2s}\rceil$ is isomorphic to $(\tilde G\lceil p^{n-s}\rceil)^{\text{t}}[p^{n-2s}]$ and thus with $H^{\text{t}}[p^{n-2s}]$.  Therefore $\{( G\lceil p^{n-s}\rceil)^{\text{t}}\lceil p^{n-2s}\rceil\}^{\text{t}}$ is isomorphic to $H[p^{n-2s}]$. From this and the fact that we have a short exact sequence
$0\rightarrow G\lceil p^{s}\rceil\rightarrow G\lceil p^{n-s}\rceil\rightarrow \{( G\lceil p^{n-s}\rceil)^{\text{t}}\lceil p^{n-2s}\rceil\}^{\text{t}}\rightarrow 0,$
we get that the Corollary \ref{C4} holds.\endproof

\subsection{Proof of Theorem \ref{T2}}

Clearly the homomorphism $g: \mathcal G[p^{n-s_Y}] \rightarrow \mathcal H[p^{n-s_Y}]$ is unique if it exists. We proved the case $Y = \Spec V$ in Corollary \ref{C3}. 
Let $y \in Y$ be a point; let $\kappa(y)$ be its residue field and let $\Spec \mathcal{R}_y$ be its local ring. If $h[p^{n-s_Y}]$ extends to a homomorphism $g_y:\mathcal G_{\mathcal{R}_y}[p^{n-s_Y}]\rightarrow\mathcal H_{\mathcal{R}_y}[p^{n-s_Y}]$, then $h[p^{n-s_Y}]$ extends also to a homomorphism  $g_{U_y}:\mathcal G_{U_y}[p^{n-s_Y}]\rightarrow\mathcal H_{U_y}[p^{n-s_Y}]$ over an open neighborhood $U_y$ of $y$ in $Y$. It follows from our assumptions that the extension $g_y$ of $h[p^{n-s_Y}]$ exists for each point $y\in Y$ of codimension 1. Indeed, if  $\text{char}\, \kappa(y) \not= p$ the group schemes $\mathcal G_{\mathcal{R}_y}$ and $\mathcal H_{\mathcal{R}_y}$ are \'etale
and therefore an extension $g_y$ trivially exists. If $\text{char}\, \kappa(y) = p$, then Corollary \ref{C3} implies that $g_y$ exists provided we have $n > s_y$ for a suitable non-negative integer $s_y$ that depends only on the ramification index of the discrete valuation ring $\mathcal{R}_y$. As the set $\Omega_p(Y)$ of points of $Y$ of codimension $1$ 
and of characteristic $p$ is finite, we can define
$$s_Y:=\max\{s_y|y\in \Omega_p(Y)\}\in\mathbb{N}.$$ 
With this $s_Y$, there exists an extension $g_U:\mathcal G_U[p^{n-s_Y}] \rightarrow \mathcal H_U[p^{n-s_Y}]$ of $h[p^{n-s_Y}]$
over an open subscheme $U \subseteq Y$ such that $\text{codim}_Y(Y\setminus U)\geq 2$.   
As $Y$ is a normal noetherian integral scheme, the existence of an extension $g:\mathcal G[p^{n-s_Y}] \rightarrow \mathcal H[p^{n-s_Y}]$ of $g_U$ is a general fact 
which holds for every locally free coherent $\mathcal{O}_Y$-modules $\mathcal{M}$ and 
$\mathcal{N}$ and for each homomorphism
$\alpha_U: \mathcal{M}_{U} \rightarrow \mathcal{N}_{U}$ of $\mathcal{O}_U$-modules. \endproof

\section{Upper bounds on heights} 

In this section we assume that $V$ is complete and $k$ is perfect. Let $\pi\in V$ and $q_{\pi}:\mathfrak{S}\twoheadrightarrow V$ be as in Section 2. We will study three different heights of a finite flat commutative group scheme $G$ of $p$ power order over $V$.  

\begin{df}\label{D3}
  {\bf (a)} By the Barsotti--Tate height of $G$ we mean the smallest non-negative integer $h_1(G)$ such that $G$ is a closed subgroup scheme of a truncated Barsotti--Tate group over $V$ of height $h_1(G)$.     

  \smallskip {\bf (b)} By the Barsotti--Tate co-height of $G$ we mean
  the smallest non-negative integer $h_2(G)$ such that $G$ is the
  quotient of a truncated Barsotti--Tate group over $V$ of height
  $h_2(G)$.

  \smallskip {\bf (c)} By the Brueil (or generator) height of $G$ we mean the
  smallest number $h_3(G)$ of generators of the $\mathfrak{S}$-module
  $N$, where $(N,\psi):=\mathbb{B}(G)$.
\end{df}

\subsection{Simple inequalities}

We have $h_2(G)=h_1(G^{\text{t}})$ and $p^{h_3(G)}$ is the order of $G_k[p]$. If $G$ is a truncated Barsotti--Tate group, then $h_1(G)=h_2(G)=h_3(G)$ are
equal to the height of $G$. Based on these properties, it is easy to check that in general we have
\begin{displaymath}
  h_3(G)=h_3(G^{\text{t}})\le\min\{h_1(G),h_2(G)\}.
\end{displaymath}
\begin{lemma}\label{L5}
  We have $h_1(G)\le 2h_3(G)$.
\end{lemma} 
{\bf Proof.} The proof of this is similar to the proof of
\cite{VZ1}, Proposition 2 (ii) but worked out in a contravariant
way. If $h_4(G)\in\{0,\ldots,h_3(G)\}$ is the smallest number of
generators of $\text{Im}(\psi)/EN$, then as in loc. cit. we argue that
there exists a Breuil window $(Q,\phi)$ relative to
$q_{\pi}:\mathfrak{S}\twoheadrightarrow V$ which has rank
$h_3(G)+h_4(G)$ and which is equipped naturally with a surjection
$(Q,\phi)\twoheadrightarrow (N,\psi)$. More precisely, starting with
$\mathfrak{S}$-linear maps $\chi_T:T:=\mathfrak{S}^{h_3(G)}\rightarrow N$ and
$\chi_L:L:=\mathfrak{S}^{h_4(G)}\rightarrow \text{Im}(\psi)$ such that
$\chi_T$ is onto and $\text{Im}(\psi)=\text{Im}(\chi_L)+EN$, one can
take $Q:=T\oplus L$ and the surjection $\chi_T\oplus
\chi_L:Q\twoheadrightarrow N$.

The existence of the surjection $(Q,\phi)\twoheadrightarrow (N,\psi)$ implies that $G$ is a closed subgroup scheme of the $p$-divisible group of height $h_3(G)+h_4(G)$ over $V$ associated to $(Q,\phi)$. Thus $h_1(G)\le h_3(G)+h_4(G)\le 2h_3(G)$.\endproof

\medskip
If $o(G)$ is as in Section \ref{2}, then we obviously have
\begin{equation}\label{F16}
  h_3(G)\le o(G).
\end{equation}
If we have a short exact sequence 
$$0\rightarrow G_1\rightarrow G\rightarrow G_2\rightarrow 0$$
of finite flat group schemes over $V$, then as the functor $\mathbb{B}$ takes short exact sequences to short exact sequences (in the category of $\mathfrak{S}$-modules endowed with Frobenius maps), we have the subadditive inequality
\begin{equation}\label{F17}
  h_3(G)\le h_3(G_1)+h_3(G_2).
\end{equation}

\begin{prop}\label{P3}
  For every truncated Barsotti--Tate group $H$ over $V$ of height $r$
  and for each $G$ whose generic fiber is isomorphic to $H_K$, we have
$$h_3(G)\le (2s+1)r.$$
Therefore we have $\max\{h_1(G),h_2(G),h_3(G)\}\le (4s+2)r$.
\end{prop} {\bf Proof.}
Let $n$ be the level of $H$. If $n\le 2s$, then from (\ref{F16}) we get that $h_3(G)\le o(G)\le 2sr$. We now assume that $n>2s$. Then from Corollary \ref{C4} we get that $G\lceil p^{n-s}\rceil/G\lceil p^{s}\rceil$ is isomorphic to $H[p^{n-2s}]$. Therefore we have $h_3(G\lceil p^{n-s}\rceil/G\lceil p^{s}\rceil)=r$. 

As the orders of $G\lceil p^{s}\rceil$ and $G/G\lceil p^{n-s}\rceil$ are equal to $sr$, from inequalities (\ref{F16}) and (\ref{F17}) we get first that $h_3(G\lceil p^{n-s}\rceil)\le h_3(G\lceil p^{n-s}\rceil/G\lceil p^{s}\rceil)+h_3(G\lceil p^{s}\rceil)\leq r+sr=(s+1)r$ and second that $h_3(G)\le h_3(G\lceil p^{n-s}\rceil)+h_3(G/G\lceil p^{n-s}\rceil)\le (s+1)r+sr=(2s+1)r$.

The group scheme $G^{\text{t}}$ satisfies the same property as $G$ i.e., $G^{\text{t}}_K$ is isomorphic to $H^{\text{t}}_K$. Thus we also have $h_3(G^{\text{t}})\le (2s+1)r$. From this and Lemma \ref{L5} we get that $h_2(G)=h_1(G^{\text{t}})\le 2h_3(G)\le (4s+2)r$. Similarly, $h_1(G)\le 2h_3(G)\le (4s+2)r$. Thus $\max\{h_1(G),h_2(G),h_3(G)\}\le (4s+2)r$.\endproof

\bigskip\smallskip\noindent
{\bf Acknowledgments.} Both authors would like to thank Bielefeld and Binghamton Universities for good working conditions. We would also like to thank the referee for many valuable comments and suggestions. The research of the junior author was partially supported by the NSF grant DMS \#0900967.


\hbox{Adrian Vasiu,\;\;\;Email: adrian@math.binghamton.edu}
\hbox{Address: Department of Mathematical Sciences, Binghamton University,} \hbox{Binghamton,  P. O. Box 6000, New York 13902-6000, U.S.A.}

\hbox{}
\hbox{Thomas Zink,\;\;\;Email: zink@math.uni-bielefeld.de}
\hbox{Address: Fakult\"at f\"ur Mathematik, Universit\"at Bielefeld,} \hbox{P.O. Box 100 131, D-33 501 Bielefeld, Germany.}

\end{document}